\documentclass[11pt,reqno]{article}
\usepackage{amsmath,amsthm}

\input{tcilatex}

\begin{document}

\begin{center}
\textbf{\large Total Operators and Inhomogeneous Proper-Value Equations}

\textbf{{\large Jos\'{e} G. Vargas}}

PST Associates, LLC (USA)\bigskip

{138 Promontory Rd, Columbia, SC 29209}\linebreak {josegvargas@earthlink.net%
\bigskip }

To Professors Eckhard Hitzer, Zbigniew Oziewicz and Wolfgang Sproessig%
{\small \bigskip }
\end{center}

\textbf{Abstract.} \ \ K\"{a}hler's two-sided angular momentum operator, $%
K+1 $, is neither vector-valued nor bivector-valued. It is total in the
sense that it involves terms for all three dimensions. Constant idempotents
that are ``proper functions'' of $K+1$'s components are not proper functions
of $K+1.$ They rather satisfy ``inhomogeneous proper-value equations'', i.e.
of the form $(K+1)U=\mu U+\pi $, where $\pi $ is a scalar.

We consider an equation of this type with $K+1$ replaced with operators $T$
that comprise $K+1$ as a factor, but also containing factors for both space
and spacetime translations. We study the action of those $T$'s on linear
combinations of constant idempotents, so that only the algebraic (spin) part
of $K+1$ has to be considered. $\pi $ is now, in general, a non-scalar
member of a K\"{a}hler algebra. We develop the system of equations to be
satisfied by the combinations of those idempotents for which $\pi $ becomes
a scalar. We solve for its solutions with $\mu =0$, which actually also
makes $\pi =0.$

The solutions with $\mu =\pi =0$ all have three constituent parts, 36 of
them being different in the ensemble of all such solutions$.$ That set of
different constituents is structured in such a way that we might as well be
speaking of an algebraic representation of quarks. In this paper, however,
we refrain from pursuing this identification in order to emphasize the
purely mathematical nature of the argument.


\section{Introduction}

This paper is a continuation of a previous one dealing with solutions of
exterior systems \cite{V55}, where we brought attention to K\"{a}hler's
concept if total angular momentum operator, $K+1$. We then introduced the
tensor product of tangent and cotangent Clifford algebras of spacetime
(cotangent here refers to differential forms viewed as functions of $r$%
-surfaces, not as antisymmetric multilinear functions of vectors). In this
new arena, we resorted to the action of the operator $K+1$, which thus makes
our study purely algebraic. We reached idempotents of the type $\mathbf{%
\epsilon }^{\pm }\boldsymbol{I}_{ij}^{\pm }\boldsymbol{P}_{l}^{\pm }$, where 
$ij$ is $(1,2)$ or $(2,3)$ or $(3,1)$ without correlations among the
superscripts. We refer readers to \cite{V55} for the definition of $\epsilon
^{\pm }$, $\boldsymbol{I}_{ij}^{\pm }$ and $\boldsymbol{P}_{l}^{\pm }$.

There are 72 such expressions when all combinations of indices are taken
into account. Symmetries reduce their number to 48. A further reduction of
those 48 to 36 follows through the study of proper values, in an extended
sense of the term, when the operator is a total operator for rotations and
spacetime translations. The steps of this study are as follows.

In section 2, we continue the discussion of the nature of total angular
momentum in K\"{a}hler's calculus of differential forms. In section 3, we
consider the emergence of algebraic inhomogeneous proper value operation
through the action of $K+1$ on idempotents. In section 4, we go one step
further in the study on the same idempotents by acting on them with the
product of $K+1$ with the operators /or translations along the axes.

In section 5, we consider the action of a second total operators, $(K+1)d%
\boldsymbol{r}$, where $d\boldsymbol{r}$ is the translation element in 3-D
Euclidean space. It is total for the displacements like $K+1$ is total for
the rotations. In section 6, we directly address the problem of finding
linear combinations of the aforementioned parts, combinations for which we
have a proper value and co-value. In the middle of the argument, we
specialize to proper value zero.

In section 7, we make our operator still more total by including time
translations. In section 8, we raise the issue of using these purely
mathematical results to formulate a physics of particles with it.

\section{K\"{a}hler's Total Angular Momentum}

We have presented in English much of K\"{a}hler's treatment (in German) of
components of angular momentum, and then of total angular momentum,
respectively in \cite{V43} and \cite{V55}. The argument starts with his
approach to Lie differentiation as a sophistical case of partial
differentiation.

Let $U_{r}$ be the differential form 
\begin{equation}
U_{r}:=a_{i_{1}\ldots i_{r}}dx^{i_{1}\ldots i_{r}},
\end{equation}%
where the coordinates are Cartesian. Consider $\partial U/\partial \phi ,$
where $\phi $ is the azimuthal coordinate. As explained by K\"{a}hler \cite%
{Kahler60}, the Lie derivative of $U_{r}$ with respect to $\phi $ is simply
the $\partial U_{r}/\partial \phi $. But $\partial U_{r}/\partial \phi $ is
not as in the literature, as we proceed to explain. We should not only
differentiate the $a_{i_{1}\ldots i_{r}}$, but also the $dx^{i}$'s, using
that $\partial (dx^{i})/\partial \phi =d(\partial x^{i}/\partial \phi )$.

Since $\partial /\partial \phi $ is of the form $f^{l}(x)\partial /\partial
x^{l}$, the preceding consideration implies that the action of operators $%
f^{l}(x)/\partial /\partial x^{l}$ on differential forms (1) have terms
additional to those where one simply differentiates the components, which is
worth remembering since this is counterintuitive. We shall not be distracted
here with details. Suffice to say that, if the $f^{l}(x)$ are constants,
then 
\begin{equation}
\left[ f^{l}\frac{\partial }{\partial x^{l}}\right] U=f^{l}\frac{\partial
a_{i_{i}\ldots i_{r}}}{\partial x^{l}}dx^{i_{1}\ldots i_{r}}.
\end{equation}%
When they are not constant, we shall distinguish between $\left[ g^{l}\frac{%
\partial }{\partial x^{m}}\right] U$ and $g^{l}\left[ \frac{\partial U}{%
\partial x^{m}}\right] $, as they are not equal. The $\partial U/\partial
x^{m}$ is in itself given by (2) with $f^{l}=0$, except for one of them, $%
f^{m}=1.$

Another important remark is that, in the K\"{a}hler calculus, $U_{r}$ is a
function of hypersurfaces. It does not depend on the coordinates used for
its description. And $\partial U_{r}/\partial \phi $ will be the total
derivative $dU_{r}/d\phi $ if $\phi $ constitutes one of the coordinates, $%
y_{n}$, in a coordinate system $(y)$ where all the other coordinates are the 
$n-1$ independent constants of the motion not additive to $y^{n}$ in the
system,%
\begin{equation}
\frac{dx^{i}}{dy^{n}}=a^{i}(x^{1},\ldots ,x^{n}).
\end{equation}%
See \cite{Kahler61}.

A moderately long process starting with the foregoing concepts leads K\"{a}%
hler to obtain the so called components of angular momentum operator in 3-D
dimensional space, $E_{3}$, as 
\begin{equation}
\frac{\partial U}{\partial \phi ^{i}}=\chi _{i}U:=\left( x^{j}\frac{\partial
U}{\partial x^{k}}-x^{k}\frac{\partial U}{\partial x^{j}}\right) +\frac{1}{2}%
w^{i}\vee U-\frac{1}{2}U\vee w^{i},
\end{equation}%
where $(i,j,k)$ are the three cyclic permutations (identity included) of $%
(1,2,3)$, where ''$\vee $'' stands for Clifford product, and where the $%
w_{i}(=w^{i})$ are the $dx^{jk}$ $(=:dx^{j}\wedge dx^{k})$. See \cite%
{Kahler60} and \cite{Kahler62}.

Except for the unit imaginary factor, the first half on the right of (4)
coincides with the components of orbital angular momentum in quantum
mechanics (with $\hbar =1$).

In view of what we have said, the $\phi ^{i}$ are azimuthal coordinates
relative to different axes and, therefore, belonging to different coordinate
systems.

K\"{a}hler defines a total angular momentum operator, $K+1$, as the
two-sided operator whose action on $U$ is given by 
\begin{equation}
(K+1)U=J_{l}Uw^{l},
\end{equation}%
with summation over repeated indices. It is linear not in the sense of
linearity in a vector space, but in the sense that it is a sum of terms each
of which contains one and only one of the $J_{l}$, the $w^{l}$ being a basis
of differential 2-forms. Applying $(K+1)$ to (5) and performing operations,
one gets 
\begin{equation}
(K+1)^{2}U=-\sum_{l}\chi _{l}^{2}U+(K+1)U.
\end{equation}%
K\"{a}hler apparently chose the symbol $K+1$ so that one gets $K(K+1)$ when
one takes to the left the last term in (6). The minus sign in front of the
first term on the right has to do with the fact that $\chi _{i}$ does not
contain the unit imaginary as factor.

We shall be applying the operator $K+1$ to Clifford valued differential
form. It is a differential operator in the K\"{a}hler subalgebra of
scalar-valued differential forms.

\section{Inhomogeneous Proper Value Equations of $K+1$}

In order to have notational continuity with \cite{V55}, we shall use the
symbol $J_{l}$ instead of $X_{l}$. Now as then, $d\boldsymbol{x}^{l}$ will
refer to $dx^{l}\mathbf{a}_{l}$ (no sum), with $l=i,j,k$ and $\mathbf{a}_{l}=%
\boldsymbol{i},\boldsymbol{j},\boldsymbol{x}.$ $d\boldsymbol{x}^{l}$ belongs
to the commutative algebra that we defined in \cite{V55}; the $J_{l}$ and $%
K+1$ do not. The action on constant differentials -- as our idempotents will
be -- of the orbital part of angular momentum operators is zero. Thus the
action of the $J_{l}$ and $K+1$ is purely algebraic. This action is zero for
scalar and for $d\boldsymbol{x}^{ijk}(\equiv dx^{ijk}\boldsymbol{a}_{i}%
\boldsymbol{a}_{j}\boldsymbol{a}_{k})$, since $dx^{ijk}$ commutes with the
whole algebra of scalar-valued differential forms and it, therefore,
commutes with $J_{l}$ and $K+1$.

In \cite{V55}, we obtained 
\begin{equation}
J_{l}d\boldsymbol{x}^{1}=(0,dx^{3},-dx^{2})\boldsymbol{a}_{1}.
\end{equation}%
We similarly have 
\begin{equation}
J_{l}d\boldsymbol{x}^{2}=(-d\boldsymbol{x}^{3},0,d\boldsymbol{x}^{1}),\;\ \
\ \ \;J_{l}d\boldsymbol{x}^{3}=(d\boldsymbol{x}^{2},-d\boldsymbol{x}x^{1},0).
\end{equation}%
Straightforward computations yields 
\begin{equation}
J_{i}d\boldsymbol{x}^{jk}=0,
\end{equation}%
and, without summation over repeated indices, 
\begin{equation}
J_{i}d\boldsymbol{x}^{ki}=\frac{1}{2}(w^{i}d\boldsymbol{x}^{ki}-d%
\boldsymbol{x}^{ki}w^{i})=-w^{k}\boldsymbol{a}_{ki}=w^{k}\boldsymbol{a}_{ik},
\end{equation}%
and 
\begin{equation}
J_{i}d\boldsymbol{x}^{ij}=\frac{1}{2}(w^{i}w^{k}-w^{k}w^{i})\boldsymbol{a}%
_{ij}=w^{j}\boldsymbol{a}_{ij}.
\end{equation}%
In general, a repeated index $i,j$ or $k$ will not mean summation, For that,
we use the index $l$, or $m$.

For later use, it is convenient to rewrite (9)-(11) by going over the three
components of the operator, acting on the same differential form: 
\begin{equation}
J_{i}d\boldsymbol{x}^{jk}=0,\;\ \ \;J_{j}d\boldsymbol{x}^{jk}=w^{k}%
\boldsymbol{a}_{jk},\;\ \ \ J_{k}d\boldsymbol{x}^{jk}=w^{j}\boldsymbol{a}%
_{kj}.
\end{equation}

As in \cite{V55}, we define 
\begin{equation}
\boldsymbol{I}_{ij}^{\pm }=\frac{1}{2}(1\pm d\boldsymbol{x}^{ij}).
\end{equation}%
Then 
\begin{equation}
J_{i}\boldsymbol{I}_{jk}^{\pm }=0,\;\ \ \;J_{i}\boldsymbol{I}_{ki}^{\pm
}=\pm \frac{1}{2}w^{k}\boldsymbol{a}_{ik},\;\ \ J_{i}\boldsymbol{I}%
_{ij}^{\pm }=\pm \frac{1}{2}w^{j}\boldsymbol{a}_{ij},
\end{equation}%
obviously without summation over repeated indices. For fixed subscripts of $%
\boldsymbol{I}^{\pm }$ and different subscripts of $J$, we have 
\begin{equation}
J_{i}\boldsymbol{I}_{jk}=0,\;\ \ \ J_{j}\boldsymbol{I}_{jk}^{\pm }=\pm \frac{%
1}{2}w^{k}\boldsymbol{x}_{jk},\;\ \ \ J_{k}\boldsymbol{I}_{jk}^{\pm }w^{j}%
\boldsymbol{a}_{jk}.
\end{equation}%
Since $w^{j}w^{i}=w^{k}$ and $w^{i}\boldsymbol{a}_{jk}=d\boldsymbol{x}^{jk}$%
, we have 
\begin{equation}
J_{j}\boldsymbol{I}_{jk}^{\pm }=\pm \frac{1}{2}w^{j}w^{i}\boldsymbol{a}%
_{jk}=-\frac{wj}{2}+\frac{wj}{2}\pm \frac{1}{2}w^{j}d\boldsymbol{x}%
^{jk}=w^{j}(\boldsymbol{I}_{jk}^{\pm }-\frac{1}{2}),
\end{equation}%
and, in a similar way, 
\begin{equation}
J_{k}\boldsymbol{I}_{jk}^{\pm }=w^{k}(\boldsymbol{I}_{jk}^{\pm }-\frac{1}{2}%
).
\end{equation}

In order to compute $(K+1)\boldsymbol{I}_{ij}^{\pm }$, it is best to use
(15) in (5): 
\begin{equation}
(K+1)\boldsymbol{I}_{ij}^{\pm }=J_{l}\boldsymbol{I}_{ij}^{\pm }w^{l}=\pm
w^{k}\boldsymbol{a}_{ij}=\pm d\boldsymbol{x}^{ij}=2\boldsymbol{I}_{ij}^{\pm
}-1,
\end{equation}%
of which we say that it is an inhomogeneous proper values equation. We refer
to $-1$ as its co-value. On the other hand 
\begin{equation}
(K+1)d\boldsymbol{x}^{ij}=(K+1)(2\boldsymbol{I}_{ij}^{+}-1)=2(K+1)%
\boldsymbol{I}_{ij}^{+}=2d\boldsymbol{x}^{ij},
\end{equation}%
which is homogeneous.

We obtained in \cite{V55}: 
\begin{equation}
(K+1)d\boldsymbol{x}^{l}=2d\boldsymbol{x}^{l},
\end{equation}%
and defined 
\begin{equation}
\boldsymbol{P}_{l}^{\pm }=\frac{1}{2}(1\pm d\boldsymbol{x}^{l}).
\end{equation}%
Clearly 
\begin{equation}
(K+1)\boldsymbol{P}_{l}^{\pm }=(K+1)(\pm \frac{d\boldsymbol{x}^{l}}{2})=\pm d%
\boldsymbol{x}^{l}=2\boldsymbol{P}_{l}^{\pm }-1,
\end{equation}%
where value and co-value are as in (18)-(19).

The $\boldsymbol{P}_{l}^{\pm }$ are proper functions of $d\boldsymbol{x}^{l}$
acting on the left. But $d\boldsymbol{x}^{l}$ is not a component of $d%
\boldsymbol{r}$; $dx^{l}$ is. The role of $d\boldsymbol{x}^{i}\boldsymbol{a}%
_{i}$ (no sum) is played by $J_{i}w^{i}$: 
\begin{equation*}
d\boldsymbol{x}^{i}U\rightarrow J_{i}Uw^{i}\qquad (\text{no sum}).
\end{equation*}

\section{First Results for Combined Rotations and Translations}

We have seen that the $\boldsymbol{I}_{ij}^{\pm }$ are proper functions with
proper value zero of $J_{k}$. But they are not proper functions of $J_{i}$, $%
J_{j}$ and $K+1$. We now study the action of $K+1$ on idempotents $%
\boldsymbol{I}_{ij}\boldsymbol{P}_{l}$. It makes a great difference whether $%
l$ is or is not equal to one of the subscripts of $\boldsymbol{I}_{ij}$. Key
for fluid computations are equations (19) and (20), as well as $(K+1)1=0$
and $(K+1)d\boldsymbol{x}^{ijk}=0.$

We easily have 
\begin{equation}
(K+1)\boldsymbol{I}_{ij}^{+}\boldsymbol{P}_{i}^{\pm }=(K+1)\frac{1}{4}%
(1+d\ldots \pm d\ldots \pm d\ldots )=\frac{1}{2}(d\ldots \pm d\ldots \pm
d\ldots ).
\end{equation}%
We have used suspension marks to indicate the irrelevance of the details
provided that none of those suspension marks is for $d\boldsymbol{x}^{ijk}$.
We thus have, adding and subtracting $1/2$,%
\begin{equation}
(K+1)\boldsymbol{I}_{ij}^{+}\boldsymbol{P}_{i}^{\pm }=\frac{1}{2}(1+d\ldots
\pm d\ldots \pm d\ldots )-\frac{1}{2}=2\boldsymbol{I}_{ij}^{\pm }%
\boldsymbol{P}_{i}^{\pm }-\frac{1}{2}.
\end{equation}%
Similarly, 
\begin{equation}
(K+1)\boldsymbol{I}_{ij}^{-}\boldsymbol{P}_{i}^{\pm }=2\boldsymbol{I}%
_{ij}^{-}\boldsymbol{P}_{i}^{\pm }-\frac{1}{2}.
\end{equation}

Consider next $(K+1)\boldsymbol{I}_{ij}^{+}\boldsymbol{P}_{k}^{\pm }.$ We
have 
\begin{align}
& (K+1)\boldsymbol{I}_{ij}^{+}\boldsymbol{P}_{k}^{\pm }=(K+1)\frac{1}{4}%
(1+d\ldots \pm d\ldots \pm d\boldsymbol{x}^{ijk})=\frac{1}{2}(d\ldots \pm
d\ldots )=  \notag \\
& =\frac{1}{2}(1+d\ldots \pm d\ldots \pm d\boldsymbol{x}^{ijk})-\frac{1}{2}%
(1\pm d\boldsymbol{x}^{ijk})=2\boldsymbol{I}_{ij}^{+}P_{k}^{\pm }-\frac{1}{2}%
(I\pm d\boldsymbol{x}^{ijk}).
\end{align}%
Proceeding similarly, we get 
\begin{equation}
(K+1)\boldsymbol{I}_{ij}^{-}\boldsymbol{P}_{k}^{\pm }=2\boldsymbol{I}%
_{ij}^{-}\boldsymbol{P}_{k}^{\pm }-\frac{1}{2}(1\pm d\boldsymbol{x}_{ijk}).
\end{equation}

We next want to know the joint effect of $K+1$ and one of the components of $%
d\boldsymbol{r}$ (in the sense of $d\boldsymbol{x}^{l}$, not $dx^{l})$.
Since $d\boldsymbol{x}^{i}\boldsymbol{I}_{ij}\boldsymbol{P}_{k}$ does
contain neither 1 nor $d\boldsymbol{x}^{ijk}$, it follows that 
\begin{subequations}
\begin{equation}
\lbrack (K+1)d\boldsymbol{x}^{i}]\boldsymbol{I}_{ij}^{+}\boldsymbol{P}%
_{k}^{\pm }=2d\boldsymbol{x}^{i}\boldsymbol{I}_{ij}^{+}\boldsymbol{P}%
_{k}^{\pm },
\end{equation}%
\begin{equation}
\lbrack (K+1)d\boldsymbol{x}^{i}]\boldsymbol{I}_{ij}^{-}\boldsymbol{P}%
_{k}^{\pm }=2d\boldsymbol{x}^{i}\boldsymbol{I}_{ij}^{+}\boldsymbol{P}%
_{k}^{\pm }.
\end{equation}

On the other hand, $d\boldsymbol{x}^{k}\boldsymbol{I}_{ij}^{+}\boldsymbol{P}%
_{i}^{\pm }$ has a $d\boldsymbol{x}^{ijk}$ term, but not a scalar term.
Hence 
\end{subequations}
\begin{subequations}
\begin{equation}
\lbrack (K+1)d\boldsymbol{x}^{k}]\boldsymbol{I}_{ij}^{+}\boldsymbol{P}%
_{i}^{\pm }=2d\boldsymbol{x}^{k}\boldsymbol{I}_{ij}^{+}\boldsymbol{P}%
_{k}^{\pm }-\frac{1}{2}d\boldsymbol{x}^{ijk},
\end{equation}%
\begin{equation}
\lbrack (K+1)d\boldsymbol{x}^{k}]\boldsymbol{I}_{ij}^{-}\boldsymbol{P}%
_{i}^{\pm }=2d\boldsymbol{x}^{k}\boldsymbol{I}_{ij}^{-}\boldsymbol{P}%
_{k}^{\pm }-\frac{1}{2}d\boldsymbol{x}^{ijk}.
\end{equation}

We now let the subscript of the $d\boldsymbol{x}$ factor inside the bracket
coincide with the subscript of $\boldsymbol{P}^{\pm }$, so that it is
absorbed. One immediately sees whether scalar and $d\boldsymbol{x}^{ijk}$
terms are present. In the absorption process, a minus sign may appear. One
readily obtains 
\end{subequations}
\begin{subequations}
\begin{equation}
\lbrack (K+1)d\boldsymbol{x}^{i}]\boldsymbol{I}_{ij}^{+}\boldsymbol{P}%
_{i}^{\pm }=\pm (2\boldsymbol{I}_{ij}^{+}\boldsymbol{P}_{i}^{\pm }-\frac{1}{2%
}),
\end{equation}%
\begin{equation}
\lbrack (K+1)d\boldsymbol{x}^{i}]\boldsymbol{I}_{ij}^{-}\boldsymbol{P}%
_{i}^{\pm }=\pm (2\boldsymbol{I}_{ij}^{-}\boldsymbol{P}_{i}^{\pm }-\frac{1}{2%
}),
\end{equation}%
\begin{equation}
\lbrack (K+1)d\boldsymbol{x}^{k}]\boldsymbol{I}_{ij}^{+}\boldsymbol{P}%
_{k}^{\pm }=\pm \lbrack 2\boldsymbol{I}_{ij}^{+}\boldsymbol{P}_{k}^{\pm }\mp 
\frac{1}{2}(1+d\boldsymbol{x}^{ijk})],  \tag{31a}
\end{equation}%
\begin{equation}
\lbrack (K+1)d\boldsymbol{x}^{k}]\boldsymbol{I}_{ij}^{-}\boldsymbol{P}%
_{k}^{\pm }=\pm \lbrack 2\boldsymbol{I}_{ij}^{-}\boldsymbol{P}_{k}^{\pm }\mp 
\frac{1}{2}(1+d\boldsymbol{x}^{ijk})].  \tag{31b}
\end{equation}

Equations (30a) and (30b) are inhomogeneous proper value equations. All the
others are not. Equations (28) are proper value equations of $(K+1)$, not of 
$(K+1)d\boldsymbol{x}^{i}$.

We have not yet integrated $K+1$ with total translation operator.

\section{Total Operator for 3-D Euclidean Space}

\setcounter{equation}{31} We shall now deal with the action of operator $%
(K+1)d\boldsymbol{r}$ on the idempotents. In principle (it would be a
tedious check), they do not commute. The question then could be: Why $(K+1)d%
\boldsymbol{r}$ instead of $d\boldsymbol{r}(K+1)$. We leave the study of the
alternative to interested readers.

For any two different indices $i$ and $j$, we have 
\end{subequations}
\begin{equation}
\boldsymbol{I}_{ij}^{+}\boldsymbol{P}_{i}^{\pm }=\boldsymbol{I}_{ij}^{+}%
\boldsymbol{P}_{j}^{\pm },\;\;\boldsymbol{OI}_{ij}^{-}\boldsymbol{P}_{i}^{+}=%
\boldsymbol{I}_{ij}^{-}\boldsymbol{P}_{j}^{\mp }.
\end{equation}%
This leads us to consider the pairs 
\begin{equation}
(\boldsymbol{I}_{12}^{+}\boldsymbol{P}_{1}^{+},\boldsymbol{I}_{12}^{+}%
\boldsymbol{P}_{2}^{+}),(\boldsymbol{I}_{12}^{+}\boldsymbol{P}_{1}^{-},%
\boldsymbol{I}_{12}^{+}\boldsymbol{P}_{2}^{-}),(\boldsymbol{I}_{12}^{-}%
\boldsymbol{P}_{1}^{+},\boldsymbol{I}_{12}^{-}\boldsymbol{P}_{2}^{-}),(%
\boldsymbol{I}_{12}^{-}\boldsymbol{P}_{1}^{-},\boldsymbol{I}_{12}^{-}%
\boldsymbol{P}_{2}^{+}),
\end{equation}%
of equal idempotents. We nevertheless do not discard half of them (which
ones in the first place?) since they are associated with different
directions for space translations. This implies that they would be
associated with different exponential factors in solutions with symmetry of
exterior systems \cite{V55}. Both copies will be used when we shall later
invoke time translations.

When doing purely algebraic computations with them, we may use, however,
just a set of four different ones. It does not matter which set we choose.
We shall consider the four $\boldsymbol{I}_{12}^{\pm }\boldsymbol{P}%
_{1}^{\ast }$ together with $\boldsymbol{I}_{12}^{\pm }\boldsymbol{P}%
_{3}^{\ast }$. We denote them as $X_{A}$, $A=1,\ldots ,8$, the ordering
being as in table 1.

The total translation operator is $d\boldsymbol{r}$, which we can write as $d%
\boldsymbol{r}^{\prime }+d\boldsymbol{x}^{3}$, where $d\boldsymbol{r}%
^{\prime }=d\boldsymbol{x}^{1}+d\boldsymbol{x}^{2}$. Notice that 
\begin{equation}
d\boldsymbol{r}^{\prime }\boldsymbol{I}_{12}^{-}=(d\boldsymbol{x}+d%
\boldsymbol{y})\boldsymbol{I}_{12}^{-}=2d\boldsymbol{x}\boldsymbol{I}%
_{12}^{+}\boldsymbol{I}_{12}^{-}=0
\end{equation}%
and also that 
\begin{equation}
d\boldsymbol{r}^{\prime }\boldsymbol{I}_{12}^{+}=d\boldsymbol{x}^{1}%
\boldsymbol{I}_{12}^{+}+d\boldsymbol{x}^{2}\boldsymbol{I}_{12}^{+}=2d%
\boldsymbol{x}^{1}\boldsymbol{I}_{12}^{+}.
\end{equation}%
Hence, 
\begin{equation}
d\boldsymbol{r}^{\prime }\boldsymbol{I}_{12}^{+}\boldsymbol{P}_{1}^{\pm }=2d%
\boldsymbol{x}^{1}\boldsymbol{I}_{12}^{+}\boldsymbol{P}_{1}^{\pm }=\pm 2%
\boldsymbol{I}_{12}^{+}\boldsymbol{P}_{1}^{\pm }.
\end{equation}%
On the basis of those considerations we build table 1. \newpage

\begin{center}
Table 1. Action of $d\boldsymbol{r}$ on the $\boldsymbol{I}_{12}^{\pm }%
\boldsymbol{P}_{1}^{\ast }$ and $\boldsymbol{I}_{12}^{\pm }\boldsymbol{P}%
_{3}^{\ast }$.\\[0pt]
\renewcommand{\arraystretch}{1.3} 
\begin{tabular}{|c|l|}
\hline
$X_{A}$ & $\ \ \ \ \ \ \ \ \ \ \ \ \ d\boldsymbol{r}^{1}X_{A}+d\boldsymbol{x}%
^{3}X_{A}$ \\ \hline
$\boldsymbol{I}_{12}^{+}\boldsymbol{P}_{1}^{+}=\frac{1}{4}(1+d\boldsymbol{x}%
^{1}+d\boldsymbol{x}^{2}+d\boldsymbol{x}^{12})$ & $\frac{1}{2}(1+d%
\boldsymbol{x}^{1}+d\boldsymbol{x}^{2}+d\boldsymbol{x}^{12})$ \smallskip \\ 
& \qquad $+\frac{1}{4}(d\boldsymbol{x}^{3}+d\boldsymbol{x}^{13}+d%
\boldsymbol{x}^{23}+d\boldsymbol{x}^{123})$ \smallskip \\ \hline
$\boldsymbol{I}_{12}^{+}\boldsymbol{P}_{1}^{-}=\frac{1}{4}(1-d\boldsymbol{x}%
^{1}-d\boldsymbol{x}^{2}+d\boldsymbol{x}^{12})$ & $\frac{1}{2}(-1+d%
\boldsymbol{x}^{1}+d\boldsymbol{x}^{2}-d\boldsymbol{x}^{12})$ \\ 
& \qquad $+\frac{1}{4}(d\boldsymbol{x}^{3}-d\boldsymbol{x}^{13}-d%
\boldsymbol{x}^{23}+d\boldsymbol{x}^{123})$ \\ \hline
$\boldsymbol{I}_{12}^{-}\boldsymbol{P}_{1}^{+}=\frac{1}{4}(1+d\boldsymbol{x}%
^{1}-d\boldsymbol{x}^{2}-d\boldsymbol{x}^{12})$ & $O+\frac{1}{4}(d%
\boldsymbol{x}^{3}+d\boldsymbol{x}^{13}-d\boldsymbol{x}^{23}-d\boldsymbol{x}%
^{123})$ \\ \hline
$\boldsymbol{I}_{12}^{-}\boldsymbol{P}_{1}^{-}=\frac{1}{4}(1-d\boldsymbol{x}%
^{1}+d\boldsymbol{x}^{2}-d\boldsymbol{x}^{12})$ & $O+\frac{1}{4}(d%
\boldsymbol{x}^{3}-d\boldsymbol{x}^{13}+d\boldsymbol{x}^{23}-d\boldsymbol{x}%
^{123})$ \\ \hline
$\boldsymbol{I}_{12}^{+}\boldsymbol{P}_{3}^{+}=\frac{1}{4}(1+d\boldsymbol{x}%
^{3}+d\boldsymbol{x}^{12}+d\boldsymbol{x}^{123})$ & $\frac{1}{2}(d%
\boldsymbol{x}^{1}+d\boldsymbol{x}^{13}+d\boldsymbol{x}^{2}+d\boldsymbol{x}%
^{23})$ \\ 
& \qquad $+\frac{1}{4}(1+d\boldsymbol{x}^{3}+d\boldsymbol{x}^{12}+d%
\boldsymbol{x}^{123})$ \\ \hline
$\boldsymbol{I}_{12}^{+}\boldsymbol{P}_{3}^{-}=\frac{1}{4}(1-d\boldsymbol{x}%
^{3}+d\boldsymbol{x}^{12}-d\boldsymbol{x}^{123})$ & $\frac{1}{2}(d%
\boldsymbol{x}^{1}-d\boldsymbol{x}^{13}+d\boldsymbol{x}^{2}-d\boldsymbol{x}%
^{23})$ \\ 
& \qquad $+\frac{1}{4}(-1+d\boldsymbol{x}^{3}-d\boldsymbol{x}^{12}+d%
\boldsymbol{x}^{123})$ \\ \hline
$\boldsymbol{I}_{12}^{-}\boldsymbol{P}_{3}^{+}=\frac{1}{4}(1+d\boldsymbol{x}%
^{3}-d\boldsymbol{x}^{12}-d\boldsymbol{x}^{123})$ & $O+\frac{1}{4}(1+d%
\boldsymbol{x}^{3}-d\boldsymbol{x}^{12}-d\boldsymbol{x}^{123})$ \\ \hline
$\boldsymbol{I}_{12}^{-}\boldsymbol{P}_{3}^{-}=\frac{1}{4}(1-d\boldsymbol{x}%
^{3}-d\boldsymbol{x}^{12}+d\boldsymbol{x}^{123})$ & $O+\frac{1}{4}(-1+d%
\boldsymbol{x}^{3}+d\boldsymbol{x}^{12}-d\boldsymbol{x}^{123})$ \\ \hline
\end{tabular}
\end{center}

We use this table in the next section.

Consider the equation 
\begin{equation}
\lbrack (K+1)d\boldsymbol{r}]X_{A}=\mu ^{\prime }X_{A}+\pi _{A}.
\end{equation}%
Needless to say that $\mu ^{\prime }$ and $\pi _{A}$ are not fully
determined and that, to the extent that they are, $\pi _{A}$ will not in
general be a number. We shall make linear combinations that are. Start by
defining 
\begin{equation}
\mu :=-\frac{\mu ^{\prime }}{4}
\end{equation}%
so that Eq. (37) can be given the form%
\begin{equation}
\lbrack (K+1)d\boldsymbol{r}+4\mu ]X_{A}=\pi _{A}.
\end{equation}%
The factor $-1/4$ in $(38)$ has been chosen to facilitate computations.

We form linear combinations $\Sigma _{A}\lambda _{A}X_{A}$. For the action
of $(K+1)d\boldsymbol{r}$ on these combinations, we make $K+1$ act on the
second column of table 1. the latter's action on $1$ and $d\boldsymbol{r}%
^{123}$ is zero, and the action on the other elements will simply multiply
them by two. We then proceed to compute 
\begin{equation*}
\lambda _{A}[(K+1)d\boldsymbol{r}+4\mu ]X_{A}.
\end{equation*}%
and sum all that up. We arrange the resulting coefficients of the $d%
\boldsymbol{x}^{l},$ $d\boldsymbol{r}^{lh}$ and $d\boldsymbol{r}^{123}$ in
columns for easy summation. We leave the scalars for last.

\begin{center}
Table 2. $\lambda _{A}[(K+1)d\boldsymbol{r}+4\mu ]X_{A}$\\[0pt]
\begin{tabular}{|c|c|c|c|c|c|c|}
\hline
$d\boldsymbol{x}^{1}$ & $d\boldsymbol{x}^{2}$ & $d\boldsymbol{x}^{3}$ & $d%
\boldsymbol{x}^{12}$ & $d\boldsymbol{x}^{13}$ & $d\boldsymbol{x}^{23}$ & $d%
\boldsymbol{x}^{123}$ \\ \hline
$\lambda _{1}+\lambda _{1}\mu $ & $\lambda _{1}+\lambda _{1}\mu $ & $\frac{%
\lambda _{1}}{2}$ & $\lambda _{1}+\lambda _{1}\mu $ & $\frac{\lambda _{1}}{2}
$ & $\frac{\lambda _{1}}{2}$ & $\frac{\lambda _{1}}{2}$ \\ \hline
$\lambda _{2}-\lambda _{2}\mu $ & $\lambda _{2}-\lambda _{2}\mu $ & $\frac{%
\lambda _{2}}{2}$ & $-\lambda _{2}+\lambda _{2}\mu $ & $-\frac{\lambda _{2}}{%
2}$ & $-\frac{\lambda _{2}}{2}$ & $\frac{\lambda _{2}}{2}$ \\ \hline
$\lambda _{3}\mu $ & $-\lambda _{3}\mu $ & $\frac{\lambda _{3}}{2}$ & $%
-\lambda _{3}\mu $ & $\frac{\lambda _{3}}{2}$ & $-\frac{\lambda _{3}}{2}$ & $%
-\frac{\lambda _{3}}{2}$ \\ \hline
$-\lambda _{4}\mu $ & $\lambda _{4}\mu $ & $\frac{\lambda _{4}}{2}$ & $%
-\lambda _{4}\mu $ & $-\frac{\lambda _{4}}{2}$ & $\frac{\lambda _{4}}{2}$ & $%
-\frac{\lambda _{4}}{2}$ \\ \hline
$\lambda _{5}$ & $\lambda _{5}$ & $\frac{\lambda _{5}}{2}+\lambda _{5}\mu $
& $\frac{\lambda _{5}}{2}+\lambda _{5}\mu $ & $\lambda _{5}$ & $\lambda _{5}$
& $\frac{\lambda _{5}}{2}+\lambda _{5}\mu $ \\ \hline
$\lambda _{6}$ & $\lambda _{6}$ & $\frac{\lambda _{6}}{2}-\lambda _{6}\mu $
& $-\frac{\lambda _{6}}{2}+\lambda _{6}\mu $ & $-\lambda _{6}$ & $-\lambda
_{6}$ & $\frac{\lambda _{6}}{2}-\lambda _{2}\mu $ \\ \hline
$0$ & $0$ & $\frac{\lambda _{7}}{2}+\lambda _{7}\mu $ & $-\frac{\lambda _{7}%
}{2}-\lambda _{7}\mu $ & $0$ & $0$ & $-\frac{\lambda _{7}}{2}-\lambda
_{7}\mu $ \\ \hline
$0$ & $0$ & $\frac{\lambda _{8}}{2}-\lambda _{8}\mu $ & $\frac{\lambda _{8}}{%
2}-\lambda _{8}\mu $ & $0$ & $0$ & $-\frac{\lambda _{8}}{2}+\lambda _{8}\mu $
\\ \hline
\end{tabular}
\end{center}

\section{Search for Solutions of Inhomogeneous Proper Values Equations of $%
(K+1)d\boldsymbol{r}$}

Solutions of inhomogeneous proper value equations of the operator $(K+1)d%
\boldsymbol{r}$ are obtained by adding each column in table 2 and setting
the sums to zero. In what way, we shall have that a surviving scalar of the
left of 
\begin{equation}
\sum_{A}\lambda _{A}[(K+1)d\boldsymbol{r}+4\mu ]X_{A}=\sum_{A}\pi _{A},
\end{equation}%
will equal the right hand side. Hence $\sum \pi _{A}$ will be the co-value
of the equation 
\begin{equation}
\lbrack (K+1)d\boldsymbol{r}][\sum_{1}^{8}\lambda _{A}X_{A}]=\mu ^{\prime
}\sum \lambda _{A}\chi _{A}+\pi ,
\end{equation}%
where $\pi (:=\sum \pi _{A}$) is a number.

Equating to zero the coefficients of $d\boldsymbol{x}^{1}$ and $d%
\boldsymbol{x}^{2}$, we respectively have 
\begin{equation}
(\lambda _{1}+\lambda _{2})+(\lambda _{5}+\lambda _{6})+\mu \lbrack (\lambda
_{1}-\lambda _{2})+(\lambda _{3}-\lambda _{4})]=0,
\end{equation}%
\begin{equation}
(\lambda _{1}+\lambda _{2})+(\lambda _{5}+\lambda _{6})+\mu \lbrack (\lambda
_{1}-\lambda _{2})-(\lambda _{3}-\lambda _{4})]=0.
\end{equation}%
The pair of equations (21)-(22) is equivalent to the pair 
\begin{equation}
\lambda _{4}=\lambda _{3}
\end{equation}%
\begin{equation}
(\lambda _{1}+\lambda _{2})+(\lambda _{5}+\lambda _{6})+\mu (\lambda
_{1}-\lambda _{2})=0.
\end{equation}

Consider next the pair of equations for $d\boldsymbol{x}^{13}$ and $d%
\boldsymbol{x}^{23}:$ 
\begin{equation}
\frac{1}{2}[(\lambda _{1}-\lambda _{2})+(\lambda _{3}-\lambda
_{4})]+(\lambda _{5}-\lambda _{6})=0,
\end{equation}%
\begin{equation}
\frac{1}{2}[(\lambda _{1}-\lambda _{2})-(\lambda _{3}-\lambda
_{4})]+(\lambda _{5}-\lambda _{6})=0.
\end{equation}%
In view of (43), they both become%
\begin{equation}
(\lambda _{1}-\lambda _{2})+2(\lambda _{5}-\lambda _{6})=0.
\end{equation}

From the terms in $d\boldsymbol{x}^{123}$ and $d\boldsymbol{x}^{3},$ we
respectively get 
\begin{equation}
(\lambda _{1}+\lambda _{2})-(\lambda _{3}+\lambda _{4})+(\lambda
_{5}+\lambda _{6})-(\lambda _{7}+\lambda _{8})+2\mu \lbrack (\lambda
_{5}-\lambda _{6})-(\lambda _{7}-\lambda _{8})]=0,
\end{equation}%
\begin{equation}
(\lambda _{1}+\lambda _{2})+(\lambda _{3}+\lambda _{4})+(\lambda
_{5}+\lambda _{6})+(\lambda _{7}+\lambda _{8})+2\mu \lbrack (\lambda
_{5}-\lambda _{6})+(\lambda _{7}-\lambda _{8})]=0.
\end{equation}%
Adding and subtracting (48) and (49), we get 
\begin{equation}
(\lambda _{1}+\lambda _{2})+(\lambda _{5}+\lambda _{6})+2\mu (\lambda
_{5}-\lambda _{6})=0
\end{equation}%
\begin{equation}
(\lambda _{3}+\lambda _{4})+(\lambda _{7}+\lambda _{8})+2\mu (\lambda
_{7}-\lambda _{8})=0.
\end{equation}%
Finally, from the equation for $d\boldsymbol{x}^{12}$, we get 
\begin{align}
& (\lambda _{1}-\lambda _{2})+\frac{1}{2}(\lambda _{5}-\lambda _{6})-\frac{1%
}{2}(\lambda _{7}-\lambda _{8})+  \notag \\
& +\mu \lbrack (\lambda _{1}+\lambda _{2})-(\lambda _{3}+\lambda
_{4})+(\lambda _{5}+\lambda _{6})-(\lambda _{7}+\lambda _{8})]=0.
\end{align}

We proceed to solve this system of equations. From (44) and (50): 
\begin{equation}
\mu (\lambda _{1}-\lambda _{2})-2\mu (\lambda _{5}-\lambda _{6})=0.
\end{equation}%
Hence, either $\mu =0$ or 
\begin{equation}
\lambda _{1}-\lambda _{2}=2(\lambda _{5}-\lambda _{6}).
\end{equation}%
We develop the option $\mu =0$ and leave the option $\mu \neq 0$ for a
future paper

Equations (47) and (54) become 
\begin{eqnarray}
\lambda _{5}+\lambda _{6} &=&-\lambda _{1}-\lambda _{2}, \\
\text{\ }\lambda _{5}-\lambda _{6} &=&-\frac{\lambda _{1}}{2}+\frac{\lambda
_{2}}{2}.
\end{eqnarray}%
Hence 
\begin{equation}
\lambda _{5}=-\frac{3}{4}\lambda _{1}-\frac{1}{4}\lambda _{2},\text{ \ \ \ \
\ \ \ \ \ \ \ \ }\lambda _{6}=-\frac{1}{4}\lambda _{1}-\frac{3}{4}\lambda
_{2}.
\end{equation}%
From (51) and (43), we obtain%
\begin{equation}
\lambda _{7}+\lambda _{8}=-2\lambda _{3},
\end{equation}%
and from (52) and (56), we further obtain%
\begin{equation}
\lambda _{7}-\lambda _{8}=\frac{3}{2}\lambda _{1}-\frac{3}{2}\lambda _{2}.
\end{equation}%
Adding and substracting, we get 
\begin{equation}
\lambda _{7}=\frac{3}{4}\lambda _{1}-\frac{3}{4}\lambda _{2}-\lambda _{3},
\end{equation}%
\begin{equation}
\lambda _{8}=-\frac{3}{4}\lambda _{1}+\frac{3}{4}\lambda _{2}-\lambda _{3}.
\end{equation}

This solution contains two simpler, complementary ones, The first one is
given by $\lambda _{2}=\lambda _{1}=1,$ $\lambda _{3}=\lambda _{4}=0.$ Then $%
\lambda _{5}=\lambda _{6}=-1,$ $\lambda _{7}=\lambda _{8}=0.$ We see three
parts in it 
\begin{equation}
\boldsymbol{I}_{12}^{+}\boldsymbol{P}_{1}^{+},\;\boldsymbol{I}_{12}^{+}%
\boldsymbol{P}_{1}^{-}(=\boldsymbol{I}_{12}^{+}\boldsymbol{P}_{2}^{-}),\;\;-(%
\boldsymbol{I}_{12}^{+}\boldsymbol{P}_{3}^{+}\oplus \boldsymbol{I}_{12}^{+}%
\boldsymbol{P}_{3}^{-}).
\end{equation}%
If we had not made $\lambda _{1}=1$, all three expressions in (63) would be
multiplied by $\lambda _{1}.$ This is not significant for our purposes.

We have used the first parenthesis to make clear that the subscripts ``$2$''
of $\boldsymbol{P}$ also is represented. And the symbol $\oplus $ is used to
indicate that, although $\boldsymbol{I}_{12}^{+}\boldsymbol{P}_{3}^{+}$ and $%
\boldsymbol{I}_{12}^{+}\boldsymbol{P}_{3}^{-}$ could be added to yield $%
\boldsymbol{I}_{12}^{+}$, we shall not do so. Thee reason is that, when
dealing with solutions of exterior systems, they will be multiplied
necessarily by different exponential factors, one with a positive exponent
and the other with a negative one.

For the second one, we make $\lambda _{1}=\lambda _{2}=0,$ $\lambda
_{3}=\lambda _{4}=1.$ It follows that $\lambda _{5}=\lambda _{6}=0,$ $%
\lambda _{7}=\lambda _{8}=-\lambda _{3}$. In this case we have 
\begin{equation}
\boldsymbol{I}_{12}^{-}\boldsymbol{P}_{1}^{+},\;\;\boldsymbol{I}_{12}^{-}%
\boldsymbol{P}_{1}^{-}(=\boldsymbol{I}_{12}^{-}\boldsymbol{P}_{2}^{+}),\;\;-(%
\boldsymbol{I}_{12}^{-}\boldsymbol{P}_{3}^{+}\oplus \boldsymbol{I}_{12}^{-}%
\boldsymbol{P}_{3}^{-}).
\end{equation}%
Once again parentheses have been used to illustrate the presence of all
these subscripts of $\boldsymbol{P}$, which justifies are speaking of these
parts.

To complete matters consider again Eq. (40). For $\mu =0$, it reduces to 
\begin{equation}
\sum_{A}\lambda _{A}(K+1)[d\boldsymbol{r}X_{A}]=\sum_{A}\pi _{A}.
\end{equation}%
$(K+1)$ acting on the non-scalars does not yield inhomogeneous terms (they
are proper functions). And it yields zero when acting on the scalars. Hence
the left hand side of (65) is zero and so is, therefore, $\sum_{A}\pi _{A}$.
Hence the co-value also is zero for those solutions with $\mu =0$. They thus
satisfy the equation 
\begin{equation}
(K+1)d\boldsymbol{r}(\sum_{A}\lambda _{A}X_{A})=0.
\end{equation}%
We could have considered this equation as our starting point and solve for
it.

We next bring (63) and (64) together and start to develop compact notation
as follows.

\begin{center}
Table 3. Constituent $\boldsymbol{I}_{12}\boldsymbol{P}$ idempotents \\[0pt]
\begin{tabular}{|c|c|c|}
\hline
$a_{1}^{3}=\boldsymbol{I}_{12}^{+}\boldsymbol{P}_{1}^{+}$ & $a_{2}^{3} = %
\boldsymbol{I}_{12}^{+}\boldsymbol{P}_{1}^{-}$ & $a_{3}^{3} = -\boldsymbol{I}%
_{12}^{\prime\,+}$ \\ \hline
$b_{1}^{3} = \boldsymbol{I}_{12}^{-}\boldsymbol{P}_{2}^{+}$ & $b_{2}^{3} = %
\boldsymbol{I}_{12}^{-}\boldsymbol{P}_{2}^{-}$ & $b_{3}^{3} = -\boldsymbol{I}%
_{12}^{\prime\, -}$ \\ \hline
\end{tabular}
\end{center}

The superscripts of $a$ and $b$ are the missing index of $I_{12}$. We used
primes on the right side of the equations for $a_{3}^{3}$ and $b_{3}^{3}$ to
keep track of the remark made between Eqs. (63) and (64) about the $%
I_{12}^{\prime \,\pm }$ representing simplifications when the exponential
factors multiplying the idempotents in solutions with symmetry of exterior
systems are neglected.

One proceeds similarly for $\boldsymbol{I}_{31}$ and $\boldsymbol{I}_{23}$.
We would thus prolong table 3 by placing table 4 above it.

\begin{center}
Table 4. Constituent $\boldsymbol{I}_{22}^{+}\boldsymbol{P}$ and $%
\boldsymbol{I}_{31}^{+}\boldsymbol{P}$ idempotents\\[0pt]
\begin{tabular}{|c|c|c|}
\hline
$a_{1}^{1}=-\boldsymbol{I}_{23}^{\prime \,+}$ & $a_{2}^{1}=\boldsymbol{I}%
_{23}^{+}\boldsymbol{P}_{2}^{+}$ & $a_{3}^{1}=\boldsymbol{I}_{23}^{+}%
\boldsymbol{P}_{2}^{-}$ \\ \hline
$b_{1}^{1}=-\boldsymbol{I}_{23}^{\prime \,-}$ & $b_{2}^{1}=\boldsymbol{I}%
_{23}^{-}\boldsymbol{P}_{3}^{+}$ & $b_{3}^{1}=\boldsymbol{I}_{23}^{-}%
\boldsymbol{P}_{3}^{-}$ \\ \hline
$a_{1}^{2}=\boldsymbol{I}_{31}^{+}\boldsymbol{P}_{3}^{-}$ & $a_{2}^{2}=-%
\boldsymbol{I}_{31}^{\prime \,+}$ & $a_{3}^{2}=\boldsymbol{I}_{31}^{+}%
\boldsymbol{P}_{3}^{+}$ \\ \hline
$b_{1}^{2}=\boldsymbol{I}_{31}^{-}\boldsymbol{P}_{1}^{-}$ & $b_{2}^{2}=-%
\boldsymbol{I}_{31}^{\prime \,-}$ & $b_{3}^{2}=\boldsymbol{I}_{31}^{-}%
\boldsymbol{P}_{1}^{+}$ \\ \hline
\end{tabular}
\end{center}

\section{Inclusion of time translations}

The operator associated with time translations is $-dt$ since the idempotent
is 
\begin{equation}
\mathbf{\epsilon }^{\pm }=\frac{1}{2}(1\mp d\boldsymbol{t})
\end{equation}%
and we have: 
\begin{equation}
-d\boldsymbol{t}\mathbf{\epsilon }^{\pm }=-d\boldsymbol{t}(1\mp d%
\boldsymbol{t})=\pm \mathbf{\epsilon }^{\pm }
\end{equation}%
The total operator now is 
\begin{equation}
T:=(-d\boldsymbol{t})(K+1)d\boldsymbol{r},
\end{equation}%
and its associated idempotents are of the type $\mathbf{\epsilon }%
\boldsymbol{I}\boldsymbol{P}$.

The idempotents given in tables 3 and 4 are now to be multiplied by $\mathbf{%
\epsilon }^{+}$ and $\mathbf{\epsilon }^{-}$. But, in the final result, we
should make equal use of the two versions of the some idempotents that enter
equations (32). The reason is that different versions of the same
idempotents will be factors in different solutions of any given system of
differential equations, where they would be multiplied by necessarily
different exponentials (called phase factors in physics). In order to
achieve that, we proceed as follows.

We reverse the signs in both superscripts of the idempotents $\boldsymbol{I}%
\boldsymbol{P}$ in table 3. We thus obtain $\boldsymbol{I}_{12}^{-}%
\boldsymbol{P}_{1}^{-},\boldsymbol{I}_{12}^{-}\boldsymbol{P}_{1}^{+},%
\boldsymbol{I}_{12}^{+}\boldsymbol{P}_{2}^{-}$ and $\boldsymbol{I}_{12}^{+}%
\boldsymbol{P}_{2}^{+}$. Rewrite these idempotents in their alternative
form: $\boldsymbol{I}_{12}^{-}\boldsymbol{P}_{2}^{+}$, $\boldsymbol{I}%
_{12}^{-}\boldsymbol{P}_{2}^{-},\boldsymbol{I}_{12}^{+}\boldsymbol{P}%
_{1}^{-} $ and $\boldsymbol{I}_{12}^{+}\boldsymbol{P}_{1}^{+}$. These
idempotents are the same ones as those that we started with, i.e. those of
table 3 though in a different order. This suggests that we multiply the
idempotents in table 3 by $\mathbf{\epsilon }^{+}$. At the same time, we
first reverse the signs in the superscripts of the original idempotents and
then multiply them by $\mathbf{\epsilon }^{-}$.

We rewrite them in a different way so that we can rewrite all of them easily
from the top of our heads. Define 
\begin{equation}
u_{l}^{m}=\mathbf{\epsilon }^{+}a_{l}^{m},\;\;\;d_{l}^{m}=\mathbf{\epsilon }%
^{+}b_{l}^{m}.
\end{equation}%
We further define $\overline{a}_{l}^{m}$ as $a_{l}^{m}$ with reversion of
the sign of each superscript. We first introduce symbols $\overline{u}%
_{l}^{m}$ and $\overline{d}_{l}^{m}$ as follows 
\begin{equation}
\overline{u}_{l}^{m}=\mathbf{\epsilon }^{-}\overline{a}_{l}^{m},\;\;%
\overline{d}_{l}^{m}=\mathbf{\epsilon }^{-}\overline{b}_{l}^{m}.
\end{equation}%
Corresponding to table 3, we would then have table 5.

\begin{center}
Table 5. Constituent idempotents of type $\epsilon \boldsymbol{I}_{1,2}%
\boldsymbol{P}$ \\[0pt]
\begin{tabular}{|c|c|c|c|}
\hline
$u/d$ & Subscript 1 & Subscript 2 & Subscript 3 \\ 
$u^{3}$ & $\boldsymbol{\epsilon}^{+}\boldsymbol{I}_{12}^{+}\boldsymbol{P}%
_{1}^{+}$ & $\boldsymbol{\epsilon}^{+}\boldsymbol{I}_{12}^{+}\boldsymbol{P}%
_{1}^{-}$ & $-\boldsymbol{\epsilon}^{+}\boldsymbol{I}_{12}^{\prime\, +}$ \\ 
\hline
$d^{3}$ & $\boldsymbol{\epsilon}^{+}\boldsymbol{I}_{12}^{-}\boldsymbol{P}%
_{2}^{+}$ & $\boldsymbol{\epsilon}^{+}\boldsymbol{I}_{12}^{-}\boldsymbol{P}%
_{2}^{-}$ & $-\boldsymbol{\epsilon}^{+}\boldsymbol{I}_{12}^{\prime\, -}$ \\ 
\hline
$\overline{d}^{3}$ & $\boldsymbol{\epsilon}^{-}\boldsymbol{I}_{12}^{+}%
\boldsymbol{P}_{2}^{-}$ & $\boldsymbol{\epsilon}^{+}\boldsymbol{I}_{12}^{+}%
\boldsymbol{P}_{2}^{+}$ & $-\boldsymbol{\epsilon}^{-}\boldsymbol{I}%
_{12}^{\prime\, +}$ \\ \hline
$\overline{u}^{3}$ & $\boldsymbol{\epsilon}^{-}\boldsymbol{I}_{12}^{-}%
\boldsymbol{P}_{1}^{-}$ & $\boldsymbol{\epsilon}^{-}\boldsymbol{I}_{12}^{-}%
\boldsymbol{P}_{1}^{+}$ & $-\boldsymbol{\epsilon}^{-}\boldsymbol{I}%
_{12}^{\prime\, -}$ \\ \hline
\end{tabular}
\end{center}

All those idempotents are different. They are not so if we remove the
factors $\boldsymbol{\epsilon}^{\pm}.$ Notice that, once we remember the top
line of this table, it is easy to reproduce the full table. We would proceed
similarly with the generation of a table of constituent idempotents of types 
$\boldsymbol{\epsilon}\boldsymbol{I}_{31}\boldsymbol{P},$ $%
\boldsymbol{\epsilon}\boldsymbol{I}_{23}\boldsymbol{P}$, helping oneself
with table 4 is needed.

\section{Concluding Remarks}

This paper is built on a few implicit premises. One of them is that there is
much to be gained from looking at solutions with symmetry of exterior
systems. Another one is that angular momentum should be preferred over
individual components. Implicit is the premise that, when going beyond total
angular momentum, as we have done, we should have $K+1$, rather than its
components, multiply other operators.

The last two of those premises would have to be abandoned to connect with
the physics as we know it. Operators $d\boldsymbol{t}/3$ and $d\boldsymbol{x}%
^{12}$ acting by product on the entries of table 5 yield familiar proper
values and so does, therefore, their sum. But this involves returning to
components. At least for spatial translations, we should keep them together.
We have kept $d\boldsymbol{t}$ separately from $d\mathbf{r}$ because boosts
does not play a role; it then appears that one should not treat $d%
\boldsymbol{t}$ and $d\mathbf{r}$ on an equal footing. It remains to be
explored what a physics based on the contents of this paper would look like.

\section{Acknowledgements}

Funding from PST\ Associates is deeply appreciated.

\end{document}